\documentclass[12pt]{article}
\newcommand{\nocopyright}{
No Copyright\thanks{
The authors hereby waive all copyright
and related or neighboring rights to this work,
and dedicate it to the public domain.
This applies worldwide.
}}
\title{Stackable and queueable permutations}
\author{Peter G. Doyle}
\date{Version 1.0 dated 30 January 2012
\\ \nocopyright
}
\newcommand{\stackit}{\mathrm{stackit}}
\newcommand{\queueit}{\mathrm{queueit}}
\usepackage{verbatim}
\usepackage{graphicx}
\usepackage{hyperref}
\newcommand{\putfig}[3]{
\begin{figure}[\floatplacement]
\centerline{\mbox{\includegraphics{figures/#2.eps}}}
\caption{#3}
\label{fig:#1}
\end{figure}
}
\begin{document}

\maketitle

\begin{abstract}
There is a natural bijection between permutations
obtainable using a stack (those avoiding the pattern $312$)
and permutations obtainable using a queue
(those avoiding $321$).
This bijection is equivalent to one described by Simion and Schmidt in 1985.
We argue that this bijection might well have been found back in 1968
by readers of volume 1
of Knuth's \emph{The Art of Computer Programming},
if Knuth had not assigned difficulty ratings to his exercises.
\end{abstract}

\section{Warm-up}

Let's warm up by playing a solitaire game called
\emph{Double-Ended Knuth};
we will call it \emph{DEK} for short,
pronounced `deek' so as to be indistinguishable from `deque' and `Dyck'.
DEK is a bare-bones relative of familiar solitaire games like Klondike.
In DEK, we use a one-suit deck consisting of only the thirteen hearts (say).
We shuffle the deck thoroughly, and place the deck face down on the table.
The goal is to end with the cards in a pile face up, running in order from 
ace to king.
In addition to \emph{the deck} and \emph{the pile} (intially empty),
we maintain a line of cards (initially empty), called \emph{the deque},
spread out face up on the board.
At any point,
if the next card needed for the pile
is available as the top card of the deck
or at either end of the deque,
we may move it up to the pile;
otherwise, our only option is to move
the top card of the deck to either end of the deque.

{\bf Exercise.}
Play enough games of DEK to decide if you like it.
This game is definitely winnable, but if you are
impatient, you might want to get rid of the face cards and play with a 10-card
deck.

\section{Introduction}
In volume 1 of \emph{The Art of Computer Programming}
\cite{knuth:vol1first},
Knuth shows that the number of permutations of $(1,2,\ldots,n)$ that
can be obtained using a stack
is the Catalan number $C_n = \frac{1}{2n+1}{{2n+1}\choose{n+1}}$.
In an exercise, he identifies the stackable permutations as those permutations
$(p_1,p_2,\ldots,p_n)$
for which there are no indices $i<j<k$ such that $p_j<p_k<p_i$.
Nowadays it is common to describe these permutations as \emph{$312$-avoiding}.

In volume 3 of \emph{The Art of Computer Programming}
\cite{knuth:vol3},
Knuth shows that the number 
of $321$-avoiding permutations
(those that contain no decreasing subsequence of length three)
is also $C_n$.
Knuth's argument here is roundabout:
Starting with a $321$-avoiding permutation,
you carry out the Robinson-Schensted algorithm
to get a pair of two-line Young tableaux of the same shape,
then fit the pair together to make a single Young tableau
of shape $(n,n)$.
Any tableau of shape $(n,n)$ arises uniquely in this way,
and the number of such tableaux is $C_n$.

There is a more straight-forward argument,
yielding a natural bijection
between $312$-avoiding and $321$-avoiding permutations.
This bijection
is equivalent to one found in 1985
by Simion and Schmidt
\cite{simionschmidt}.
(See Claesson and Kitaev \cite{claesson}.)
We argue that this bijection might well have been found 
back in 1968
by readers of Knuth's volume 1,
if Knuth had not taken care to grade all his exercises
as to the expected difficulty of solving them.

\section{One example}

Maybe one example will allow you to skip the rest of this.
Here are two corresponding permutations,
the first obtained using a stack and the second obtained using a queue.
Figure \ref{fig:weakdyck} shows the corresponding stack height and queue
length record.
\putfig{weakdyck}{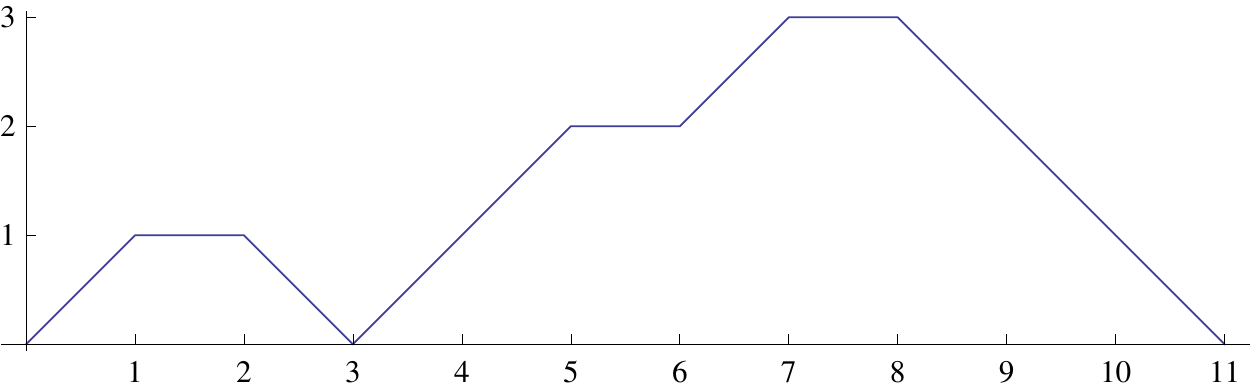}{Stack height and queue length:
A peakless weak Dyck path}
\[
\begin{array}{l|r|r}
\mbox{input}
&\mbox{stack}
&\mbox{output}
\\
\hline

& 

& 
1, 2, 3, 4, 5, 6, 7
\\

& 
1
&
2, 3, 4, 5, 6, 7
\\
2
&
1
&
3, 4, 5, 6, 7
\\
2, 1
& 

&
3, 4, 5, 6, 7
\\
2, 1
&
3
&
4, 5, 6, 7
\\
2, 1
&
4, 3
&
5, 6, 7
\\
2, 1, 5
&
4, 3
&
6, 7
\\
2, 1, 5
&
6, 4, 3
&
7
\\
2, 1, 5, 7
&
6, 4, 3
&

\\
2, 1, 5, 7, 6
&
4, 3
&

\\
2, 1, 5, 7, 6, 4
&
3
&

\\
2, 1, 5, 7, 6, 4, 3
&

&

\end{array}
\]
\[
\begin{array}{l|r|r}
\mbox{output}
&
\mbox{queue}
&
\mbox{input}
\\
\hline
 
&

&
1, 2, 3, 4, 5, 6, 7 
\\
 
&
1
&
2, 3, 4, 5, 6, 7
\\
2
&
1
&
3, 4, 5, 6, 7
\\
2, 1
& 

&
3, 4, 5, 6, 7
\\
2, 1
&
3
&
4, 5, 6, 7
\\
2, 1
&
3, 4
&
5, 6, 7
\\
2, 1, 5
&
3, 4
&
6, 7
\\
2, 1, 5
&
3, 4, 6
&
7
\\
2, 1, 5, 7
&
3, 4, 6
&

\\
2, 1, 5, 7, 3
&
4, 6
&

\\
2, 1, 5, 7, 3, 4
&
6
&

\\
2, 1, 5, 7, 3, 4, 6
& 

&

\\
\end{array}
\]

\section{An exercise rated ``M28''}

Toward the beginning of volume 1 of \emph{The Art of Computer Programming}
\cite{knuth:vol1first},
Knuth poses the following exercise.
The rating ``M28'' identifies this exercise as
mathematically oriented,
and of average to moderate difficulty, likely requiring
somewhere between fifteen to twenty minutes (difficulty ``20'') and 
over two hours (difficulty ``30'') to solve.
(The scale is `logarithmic'.)

{\bf Exercise 2.2.1-5 [M28]}
\emph{
Show that it is possible to obtain the permutation
$p_1\,p_2\,\ldots\,p_n$
from $1\,2\,\ldots\,n$
using a stack
if and only if there are no indices $i<j<k$ such that $p_j<p_k<p_i$.
}

{\bf Solution.}
This exercise refers to previous exercises,
which make clear that we are to produce the permutation
by a sequence of push and pop operations,
where the objects pushed are the numbers
$1,2,\ldots,n$, in that order,
and the output is the sequence of objects popped.
For example,
by pushing and popping in strict alternation,
we get the identity permutation $(1,2,\ldots,n)$.
By doing $n$ pushes followed by $n$ pops,
we get the reversing permutation $(n,\ldots,2,1)$.
Let us call any permutation obtainable in this way \emph{stackable}.

Rephrased in the language of pattern-avoiding permutations,
the problem is to show that a permutation is stackable
if and only if it avoids the pattern $312$.

This condition is necessary
because to obtain a permutation containing a pattern of form $312$,
when we pop $3$, we would have to have $1$ and $2$ on the stack,
with $1$ above $2$, which is impossible.

To see that the condition is sufficient,
suppose that we attempt to produce a given output permutation
in the obvious way---the only way possible---namely,
working through the desired output permutation in order,
pushing to get the next desired output onto the top of the stack
if it is not already on the stack,
and hoping that it is at the top of the stack
if it is already on the stack.
The only thing that can go wrong
is if at some point when we wish to take something
(call it $1$) from the stack,
there is something larger (call it $2$) blocking it;
if $1$ and $2$ are on the stack,
it is because we have already popped something larger still (call it $3$);
so $3$ has been called first, now $1$ is called, with $2$ to be called later.

If we associate to a stackable permutation
the record of stack height as a function of time,
we get a Dyck path of length $2n$,
meaning a path in the integer lattice $\mathbf{Z}^2$ from $(0,0)$ to $(2n,0)$,
made up of steps $(1,1)$ and $(1,-1)$, that never goes below the $x$-axis.
Any such Dyck path arises
from one and only one obtainable permutation.
The number of Dyck paths of length $2n$ is the Catalan number
$C_n = \frac{1}{2n+1}{{2n+1}\choose{n+1}}$.
(Pick a sequence consisting of $n+1$ $1$'s and $n$ $-1$'s;
take the unique rotation for which all partial sums are non-negative;
discard the initial $1$; convert to a Dyck path.)
Thus the number of stackable permutations
is $C_n$.

Notice that the answer to this problem does not change if,
in addition to pushing and popping,
we have the option of transferring an object
directly from the input to the output.
Now there is more than one sequence of operations
that produces a given permutation,
because we have the option of transferring an object
directly from input to output,
or pushing it and then immediately popping it.
To avoid this ambiguity,
we may assume that we never pop immediately after pushing,
say because we wish to minimize the number of operations.
Now in the stack record of an obtainable permutation,
all of the peaks
(consisting of a step $(1,1)$ followed immediately by a step $(1,-1)$)
will have been replaced by horizontal steps $(1,0)$;
the total number of steps in the path will now be shorter than $2n$
by the number of peaks removed.
The $C_n$ possible stack records will now consist of
`peakless weak Dyck paths'.

\section{An exercise rated ``00''}

Immediately following the exercise we've just discussed
comes this throw-away exercise, with difficulty rating ``00''.

{\bf Exercise 2.2.1-6 [00]}
\emph{
Consider the problem of exercise 2, with a queue substituted for a stack.
What permutations of $1\,2\,\ldots\,n$ can be obtained with the use of a queue?
}

{\bf Solution.}
The set-up of exercise 2 is the same as that of exercise 5, the preceeding
exercise,
which we've just discussed.
If we interpret the problem ungenerously,
as Knuth's ``00'' rating indicates that he expects us to do,
the only permutation that can be obtained is the identity permutation.

Let's vary the problem by allowing direct transfers from input to output.
Surely this would be a reasonable alternative interpretation of
permutations that `can be obtained with the use of a queue'.
While allowing direct transfers does not change which permutations can
be obtained with the use of a stack, here we find that we can obtain
all and only those permutations that arise by interleaving two increasing
subsequences.
Let us call permutations that can be obtained using a queue in this way
\emph{queueable}.

A permutation is queueable just if it is $321$-avoiding.
The reason is that a permutation
can be obtained by interleaving two increasing subsequences just if
has no decreasing subsequence of length three.
This is a simple consequence of the Robinson-Schensted
correspondence, but we do not need anything that fancy.
Any permutation obtained by interleaving two increasing subsequences
must avoid $321$, because if not, one of the two subsequences would contain
the pattern $21$, making it non-increasing.
The argument in the other direction is just like the argument above that
a $312$-avoiding permutation is stackable.
Here we want to show that a $321$-avoiding permutation is queueable.
Once again, 
we work through the desired output permutation,
hoping that the next desired output is at the head of the queue
if it is already in the queue, and moving items to the queue to uncover
it if it is still in the input stream.
The only thing that can go wrong is if when we want to take something
(call it $2$) from the queue,
there is something larger (call it $3$) blocking it;
if $2$ is in the queue,
it is because we had to make way to retrieve something smaller 
(call it $1$) from the input stream.
So in the input $3$ was ahead of $2$, and $2$ was ahead of $1$.

Here, as in the case of a stack when we allowed direct transfers from
input to output,
a given queueable permutation can arise in multiple ways,
because we can move an object to the queue early or late.
To avoid this ambiguity,
let us assume that we never move an object to the queue unless and until
we have to, as in the procedure we've just described.
In contrast to what we saw with stacks,
the preferred way of obtaining a permutation is not necessarily quicker
than competing ways, since moving an object to the queue early does not
increase the number of operations.
But the preferred way does minimize the aggregate time spent by objects on the
queue, so if we imagine that a storage fee is charged for the queue,
the preferred way minimizes this fee.
Note that minimizing total time in storage works to pick out
the preferred way of using a stack, as well as a queue.

In picking out a preferred way to obtain
a queueable permutation,
we implicitly pick out a canonical way to decompose it into
(at most) two increasing subsequences.
One subsequence, consisting of those objects that are transferred directly,
consists of the `record-setting' objects, namely, those that are larger
than any predecessor.
This subsequence is of course increasing, no matter what permutation
we start with.
A permutation is queueable just if the complemetary subsequence is also
increasing.

Now recall that in the case of stackable permutations,
where we allow direct transfers,
the stack height record is a weak Dyck path;
any weak Dyck path determines a unique stackable permutation;
and there is a bijection between stackable permutations and
peakless weak Dyck paths.
The same holds here, except that now instead of stack height we look
at queue length:
The queue length record is a weak Dyck path;
any weak Dyck path determines a unique queueable permutation;
and there is a bijection between queueable permuations and
peakless weak Dyck paths,
associating to any queueable permutation
the queue length record of the preferred way of obtaining it.
The main difference is in the
way you would go about reducing a general weak Dyck path to the
unique peakless weak Dyck path realizing the same queueable permutation.

Because of the correspondence to peakless weak Dyck paths,
we get that the number of queueable permutations is $C_n$,
just as for stackable permutations.
Furthermore, we get a natural bijection between
stackable and queueable permutations by pairing up those that share
the same peakless weak Dyck path.
This bijection is equivalent to that found by
Simion and Schmidt
\cite{simionschmidt}.
Claesson and Kitaev
\cite{claesson}
call this bijection
(or something equivalent to it)
the \emph{Knuth-Richards bijection} for the following reason.
If we turn a peakless weak Dyck path into a standard Dyck path
in the obvious way,
the correspondence between queueable permutations and standard Dyck paths
that we get is equivalent to that described by Richards
\cite{richards:knuth}.
Turning a stackable permuation into a Dyck path following Knuth,
and turning the Dyck path into a queueable permutation following
Richards,
we get a bijection between stackable and queueable permutations.

In describing the bijection we have found as being mediated either by
peakless weak Dyck paths or standard Dyck paths,
we are not doing it justice.
We can describe it more directly by saying that,
if you set out to produce a given stackable permutation using a stack,
but by mistake you use a queue instead, you get the corresponding
queueable permutation, and vice versa.

Here's another way to think of it.
Suppose that we have at our disposal not a stack or a queue, but a \emph{set},
from which we are able to recover objects in any order.
Now we can realize any permutation.
Again, let us imagine that there is a cost for using the set,
so we use it as little as possible.
To any permutation there corresponds a unique peakless weak Dyck path,
namely the set size record,
only now the correspondence is many-to-one.
Start with an arbitrary permutation $\sigma$,
find the corresponding peakless weak Dyck path,
and let $\stackit(\sigma)$ be the corresponding stackable permutation.
The function $\stackit$ is an idempotent map
($\stackit \circ \stackit = \stackit$),
projecting the set of all permutations onto the set of stackable
permutations.
The mapping
$\stackit$ tells what you get if you try to produce a specified permutation
by means of a set, but unbeknownst to you the set is not a set but a stack,
and you always get the top element no matter which one you ask for.
Similarly, we get a idempotent function $\queueit$ telling what you get
if what you thought was a set was really a queue.
Restricting $\queueit$ to stackable permutations gives a bijection to
queueable permutations;
its inverse is the restriction of $\stackit$ to queueable permutations.

One more thing:
It is plausible to say that using a queue, we can obtain any permuation,
because after all, we can retrieve any object from the queue by succesively
transferring cards from the output back to the input until we find the
object we are looking for.
So here is another interpretation of the problem that would justify a
difficulty rating of ``00''.

\section{Conclusion}

The drawback of labelling exercises according to the expected difficulty is
that it does not allow for ambiguous exercises, where the level of difficulty
may depend on how the exercise is interpreted.
Assigning exercises is not like writing a computer program, where
ambiguity is to be avoided at all costs.
The problems that the real world poses for us are almost always ambiguous.
And even if unambiguous, how difficult they may be is seldom known in advance.

\section{Problems}

\begin{enumerate}
\item
What permutations can be obtained using two stacks?
\item
What permutations can be obtained using two queues?
\item
What permutations can be obtained using one stack and one queue?
\item
What is the probability of winning at DEK if you play optimally?
\nocite{evenitai}
\nocite{pratt:dek}
\nocite{tarjan:networks}
\nocite{rt:dek}
\nocite{bona:survey}
\nocite{unger:88}
\nocite{unger:92}
\end{enumerate}

\bibliography{knuth}
\bibliographystyle{hplain}

\end{document}